\documentclass[final]{siamart1116}

\usepackage{amsfonts}
\usepackage{mathrsfs}
\usepackage{graphicx}
\usepackage{epstopdf}
\usepackage[algo2e,ruled,linesnumbered,noline,algosection]{algorithm2e}
\usepackage{comment}
\usepackage{color}
\usepackage{tabto}
\usepackage{multirow}
\usepackage[caption=false]{subfig}
\ifpdf
  \DeclareGraphicsExtensions{.eps,.pdf,.png,.jpg}
\else
  \DeclareGraphicsExtensions{.eps}
\fi

\newsiamremark{remark}{Remark}

\crefname{algocf}{Algorithm}{Algorithms}
\Crefname{algocf}{Algorithm}{Algorithms}

\renewcommand{\theequation}{\thesection.\arabic{equation}}
\renewcommand{\thefigure}{\thesection.\arabic{figure}}
\renewcommand{\thetable}{\thesection.\arabic{table}}

\numberwithin{theorem}{section}
\numberwithin{equation}{section}
\numberwithin{figure}{section}
\numberwithin{table}{section}
\numberwithin{algorithm}{section}

\title{Low-Rank Solution Methods for Stochastic Eigenvalue Problems%
   \thanks{This work was supported by the U.S. Department of Energy Office of Advanced Scientific Computing Research, Applied Mathematics program under award DE-SC0009301 and by the U.S. National Science Foundation under grant DMS1418754.}}

\author{
   Howard C. Elman%
   \thanks{Department of Computer Science and Institute for Advanced Computer Studies, University of Maryland, College Park, MD 20742 (\email{elman@cs.umd.edu}).}
   \and
   Tengfei Su%
   \thanks{Applied Mathematics \& Statistics, and Scientific Computation Program, University of Maryland, College Park, MD 20742 (\email{tengfesu@math.umd.edu}).}
}

\headers{Low-Rank Solution Methods for Stochastic Eigenvalue Problems}
{H. C. Elman and T. Su}

\begin{document}

\maketitle

\begin{abstract} 
We study efficient solution methods for stochastic eigenvalue problems arising from discretization of self-adjoint partial
differential equations with random data. With the stochastic Galerkin approach, the solutions are represented as generalized polynomial chaos expansions. A low-rank variant of the inverse subspace iteration algorithm is presented for computing one or several minimal eigenvalues and corresponding eigenvectors of parameter-dependent matrices. In the algorithm, the iterates are approximated by low-rank matrices, which leads to significant cost savings.  The algorithm is tested on two benchmark problems, a stochastic diffusion problem with some poorly separated eigenvalues, and an operator derived from a discrete stochastic Stokes problem whose minimal eigenvalue is related to the inf-sup stability constant. Numerical experiments show that the low-rank algorithm produces accurate solutions compared to the Monte Carlo method, and it uses much less computational time than the original algorithm without low-rank approximation.

\end{abstract}



\section{Introduction}
\label{sec:intro}
 
Approaches for solving stochastic eigenvalue problems can be broadly divided into non-intrusive methods, including Monte Carlo methods and stochastic collocation methods \cite{AnSc12,PrSc02}, and intrusive stochastic Galerkin methods. The Galerkin approach gives parametrized descriptions of the eigenvalues and eigenvectors, represented as expansions with stochastic basis functions. A commonly used framework is the generalized polynomial chaos (gPC) expansion \cite{XiKa03}. A direct projection onto the basis functions will result in large coupled nonlinear systems that can be solved by a Newton-type algorithm \cite{BeOn17,GhGh07}. Alternatives that do not use nonlinear solvers are stochastic versions of the (inverse) power methods and subspace iteration algorithms \cite{HaKa15,HaLa17,MeGh14,SoEl16,VeGu06}. These methods have been shown to produce accurate solutions compared with the Monte Carlo or collocation methods. However, due to the extra dimensions introduced by randomness, solving the linear systems, as well as other computations, can be expensive. In this paper, we develop new efficient solution methods that use low-rank approximations for the stochastic eigenvalue problems with the stochastic Galerkin approach.

Low-rank methods have been explored for solution of stochastic/parametrized partial differential equations (PDEs) and high-dimensional PDEs. Discretization of such PDEs gives large, sparse, and in general structured linear systems. Iterative solvers construct approximate solutions of low-rank matrix or tensor structure so that the matrix-vector products can be computed cheaply. Combined with rank compression techniques, the iterates are forced to stay in low-rank format. This idea has been used with Krylov subspace methods \cite{BaGr13,BeOn15,KrTo11,LeEl17} and multigrid methods \cite{ElSu17,Hack15}. The low-rank solution methods solve the linear systems to a certain accuracy with much less computational effort and facilitate the treatment of larger problem scales. Low-rank iterative solvers were also used in \cite{BeDo16,BeDo17} for optimal control problems constrained by stochastic PDEs. 

In this study, we use the stochastic Galerkin approach to compute gPC expansions of one or more minimal eigenvalues and corresponding eigenvectors of parameter-dependent matrices, arising from discretization of stochastic self-adjoint PDEs. Our work builds on the results in \cite{MeGh14,SoEl16}. We devise a low-rank variant of the stochastic inverse subspace iteration algorithm, where the iterates and solutions are approximated by low-rank matrices. In each iteration, the linear system solves required by the inverse iteration algorithm are performed by low-rank iterative solvers. The orthonormalization and Rayleigh quotient computations in the algorithm are also computed with the low-rank representation. To test the efficiency of the proposed algorithm, we consider two benchmark problems, a stochastic diffusion problem and a Schur complement operator derived from a discrete stochastic Stokes problem. The diffusion problem has some poorly separated eigenvalues and we show that a generalization of Rayleigh-Ritz refinement for the stochastic problem can be used to obtain good approximations. A low-rank geometric multigrid method is used for solving the linear systems. For the Stokes problem, the minimal eigenvalue of the Schur complement operator is the square of the parametrized inf-sup stability constant for the Stokes operator. Each step of the inverse iteration entails solving a Stokes system for which a low-rank variant of the MINRES method is used. We demonstrate the accuracy of the solutions and efficiency of the low-rank algorithms by comparison with the Monte Carlo method and the full subspace iteration algorithm without using low-rank approximation.

We note that a low-rank variant of locally optimal block preconditioned conjugate gradient method was studied in \cite{KrTo11_2} for eigenvalue problems from discretization of high-dimensional elliptic PDEs. Another dimension reduction technique is the reduced basis method. This idea is used in \cite{FuMa16,HoWo17,MaMa00}, where the eigenvectors are approximated from a linear space spanned by carefully selected sample ``snapshot'' solutions via, for instance, a greedy algorithm that minimizes an a posteriori error estimator. Inf-sup stability problems were also studied in \cite{HuRo07,SiKr16} in which lower and upper bounds for the smallest eigenvalue of a stochastic Hermitian matrix are computed using successive constraint methods in the reduced basis context.

The rest of the paper is organized as follows. In \cref{sec:stoch_isi} we review the stochastic inverse subspace iteration algorithm for computing several minimal eigenvalues and corresponding eigenvectors of parameter-dependent matrices. In \cref{sec:lowrank} we introduce the idea of low-rank approximation in this setting, and discuss how computations in the inverse subspace iteration algorithm are done efficiently with quantities in low-rank format. The stochastic diffusion problem and the stochastic Stokes problem are discussed in  \cref{sec:diff,sec:stokes}, respectively, with numerical results showing the effectiveness of the low-rank algorithms. Conclusions are drawn in the last section.

\section{Stochastic inverse subspace iteration}
\label{sec:stoch_isi}

Let $(\Omega, \mathcal{F}, \mathcal{P})$ be a probability triplet where $\Omega$ is a sample space with $\sigma$-algebra $\mathcal{F}$ and probability measure $\mathcal{P}$. Define a random variable $\xi:\Omega\rightarrow \Gamma\subset\mathbb{R}^m$ with uncorrelated components and let $\mu$ be the induced measure on $\Gamma$. Consider the following stochastic eigenvalue problem: find $n_e$ minimal eigenvalues $\lambda^s(\xi)$ and corresponding eigenvectors $u^s(\xi)$ such that
\begin{equation} 
    \label{eq:eigpb}
    A(\xi)u^s(\xi) = \lambda^s(\xi)u^s(\xi),\quad s=1,2,\ldots,n_e,
\end{equation}
almost surely, where $A(\xi)$ is a matrix-valued random variable. We will use a version of stochastic inverse subspace iteration studied in \cite{MeGh14,SoEl16} for solution of \cref{eq:eigpb}. The approach derives from a stochastic Galerkin formulation of subspace iteration, which is based on projection onto a finite-dimensional subspace of $L^2(\Gamma)$ spanned by the gPC basis functions $\{\psi_k(\xi)\}_{k=1}^{n_\xi}$. These functions are orthonormal, with
\begin{equation}
    \langle \psi_i\psi_j\rangle = \mathbb{E}[\psi_i\psi_j] = \int_\Gamma \psi_i(\xi)\psi_j(\xi)\text{d}\mu = \delta_{ij},
\end{equation}
where $\langle\cdot\rangle$ is the expected value, and $\delta_{ij}$ is the Kronecker delta. The stochastic Galerkin solutions are expressed as expansions of the gPC basis functions,
\begin{equation}
    \label{eq:gpc_expand}
    \lambda^s(\xi)=\sum_{r=1}^{n_\xi} \lambda^s_r\psi_r(\xi), \quad {u}^s(\xi)=\sum_{j=1}^{n_\xi} {u}^s_j\psi_j(\xi). 
\end{equation}

We briefly review the stochastic subspace iteration method in the case where $A(\xi)$ admits an affine expansion with respect to components of the random variable $\xi$:
\begin{equation}
    \label{eq:A_expansion}
    A(\xi) = A_0 + \sum_{l=1}^m A_l\xi_l
\end{equation}
where each $A_l$ is an $n_x\times n_x$ deterministic matrix, obtained from, for instance, finite element discretization of a PDE operator. The matrix $A_0$ is the mean value of $A(\xi)$. Such a representation can be obtained from a Karhunen-Lo\`{e}ve (KL) expansion \cite{Loeve} of the stochastic term in the problem (see \cref{eq:kl}). Let $\{u^{s,(i)}(\xi)\}_{s=1}^{n_e}$ be a set of approximate eigenvectors obtained at the $i$th step of the inverse subspace iteration. Then at step $i+1$, one needs to solve
\begin{equation}
    \label{eq:ii_proj}
    \langle Av^{s,(i+1)}\psi_k\rangle = \langle u^{s,(i)}\psi_k\rangle, \quad k=1,2,\ldots,n_\xi,
\end{equation}
for $\{v^{s,(i+1)}\}_{s=1}^{n_e}$ and compute $\{u^{s,(i+1)}\}_{s=1}^{n_e}$ via orthonormalization. If $n_e=1$, for the latter requirement, $v^{s,(i+1)}$ is normalized so that $||u^{s,(i+1)}||_2=1$ almost surely. If $n_e>1$, a stochastic version of the Gram-Schmidt process is applied and the resulting vectors $\{u^{s,(i+1)}\}_{s=1}^{n_e}$ satisfy $\langle u^{s,(i+1)},u^{t,(i+1)}\rangle_{\mathbb{R}^{n_x}} = \delta_{st}$ almost surely, where $\langle\cdot,\cdot\rangle_{\mathbb{R}^{n_x}}$ is the Euclidean inner product in $\mathbb{R}^{n_x}$. With the iterates expressed as gPC expansions, for instance, ${u}^{s,(i)}(\xi)=\sum_{j=1}^{n_\xi} {u}^{s,(i)}_j\psi_j(\xi)$, collecting the $n_\xi$ equations in \cref{eq:ii_proj} for each $s$ yields an $n_x n_\xi\times n_x n_\xi$ linear system
\begin{equation}
    \label{eq:ii_solve}
    \sum_{l=0}^{m} (G_l\otimes A_l)\mathbf{v}^{s,(i+1)} = \mathbf{u}^{s,(i)}
\end{equation}
where $\otimes$ is the Kronecker product, each $G_l$ is an $n_\xi\times n_\xi$ matrix with $[G_l]_{kj}=\langle\xi_l\psi_k\psi_j\rangle$ ($\xi_0\equiv 1$ and $G_0=I$), and
\begin{equation}
    \mathbf{u}^{s,(i)} = \begin{pmatrix}
    {u}_1^{s,(i)}\\
    {u}_2^{s,(i)}\\
    \vdots\\
    {u}_{n_\xi}^{s,(i)}\\
    \end{pmatrix} \in \mathbb{R}^{n_xn_\xi}.
\end{equation}
Note that the matrices $\{G_l\}$ are sparse due to orthogonality of the gPC basis functions \cite{ErUl10,PoEl09}. The initial iterate is given by solving the mean problem $A_0\bar u^s=\bar\lambda^s\bar u^s$, and
\begin{equation}
    \mathbf{u}^{s,(0)} = \begin{pmatrix}
    \bar u^s\\
    0\\
    \vdots\\
    0\\
    \end{pmatrix}.
\end{equation}
The complete algorithm is summarized as \cref{alg:sisi}. The details of the computations in steps 4 and 7 are given in \cref{eq:normal_vec,eq:gram-schmidt,eq:rayleigh} below. 

\begin{algorithm2e}
\caption{Stochastic inverse subspace iteration}
\label{alg:sisi}
\SetNlSty{}{}{:}
\DontPrintSemicolon
\SetKw{Init}{initialization}
\Init: initial iterate $\mathbf{u}^{s,(0)}$. \;
\For{$i=0,1,2,\ldots$}{
    Solve the stochastic Galerkin system \cref{eq:ii_solve} for $\mathbf{v}^{s,(i+1)}$, $s=1,2,\ldots,n_e$. \;
    If $n_e=1$, compute $\mathbf{u}^{s,(i+1)}$ by normalization. Otherwise, apply a stochastic Gram-Schmidt process for orthonormalization. \;
    Check convergence.
}
Compute eigenvalues using a Rayleigh quotient.
\end{algorithm2e}

\section{Low-rank approximation}
\label{sec:lowrank}

In this section we discuss the idea of low-rank approximation and how this can be used to reduce the computational costs of \cref{alg:sisi}. The size of the Galerkin system \cref{eq:ii_solve} is in general large and solving the system can be computationally expensive. We utilize low-rank iterative solvers where the iterates are approximated by low-rank matrices and the system is efficiently solved to a specified accuracy. In addition, low-rank forms can be used to reduce the costs of the orthonormalization and Rayleigh quotient computations in the algorithm.

\subsection{System solution}
For any random vector $x(\xi)$ with expansion $x(\xi)=\sum_{j=1}^{n_\xi} x_j\psi_j(\xi)$ where each $x_j$ is a vector of length $n_x$, let
\begin{equation}
    X = \text{mat}(\mathbf{x}) = [x_1,x_2,\ldots,x_{n_\xi}] \in \mathbb{R}^{n_x\times n_\xi}.
\end{equation}
Then the Galerkin system $\sum_{l=0}^{m} (G_l\otimes A_l)\mathbf{x}= \mathbf{f}$ is equivalent to the matrix form
\begin{equation}
    \label{eq:system_matrix}
    \sum_{l=0}^{m} A_l XG_l^T = F = \text{mat}(\textbf{f}).
\end{equation}
Let $X^{(i)}=\text{mat}(\mathbf{x}^{(i)})$ be the $i$th iterate computed by an iterative solver applied to \cref{eq:system_matrix}, and suppose $X^{(i)}$ is represented as the product of two rank-$\kappa$ matrices, i.e., $X^{(i)}=Y^{(i)}Z^{(i)T}$, where $Y^{(i)}\in\mathbb{R}^{n_x\times\kappa}$, $Z^{(i)}\in\mathbb{R}^{n_\xi\times\kappa}$. If this factored form is used throughout the iteration without explicitly forming $X^{(i)}$, then the matrix-vector product $(G_l\otimes A_l)\mathbf{x}$ will have the same structure,
\begin{equation}
A_l X^{(i)}G_l^T = (A_lY^{(i)})(G_lZ^{(i)})^T,
\end{equation}
and it is only necessary to compute $A_lY^{(i)}$ and $G_lZ^{(i)}$. If $\kappa\ll\min(n_x,n_\xi)$, this means that the computational costs of the matrix operation are reduced from $O(n_xn_\xi)$ to $O((n_x+n_\xi)\kappa)$. On the other hand, summing terms with the factored form tends to increase the rank, and rank compression techniques must be used in each iteration to force the matrix rank $\kappa$ to stay low. In particular, if $X_1^{(i)}=Y_1^{(i)}Z_1^{(i)T}$, $X_2^{(i)}=Y_2^{(i)}Z_2^{(i)T}$, where $Y_1^{(i)}\in\mathbb{R}^{n_x\times\kappa_1}$, $Z_1^{(i)}\in\mathbb{R}^{n_\xi\times\kappa_1}$, $Y_2^{(i)}\in\mathbb{R}^{n_x\times\kappa_2}$, $Y_2^{(i)}\in\mathbb{R}^{n_\xi\times\kappa_2}$, then
\begin{equation}
X_1^{(i)} + X_2^{(i)} = [Y_1^{(i)},Y_2^{(i)}][Z_1^{(i)},Z_2^{(i)}]^T.
\end{equation}
The addition gives a matrix of rank $\kappa_1+\kappa_2$ in the worst case. Rank compression can be achieved by an SVD-based truncation operator $\tilde X^{(i)} = \mathcal{T}(X^{(i)})$ so the matrix $\tilde X^{(i)}$ has a much smaller rank than $X^{(i)}$ \cite{KrTo11}. Specifically, we compute QR factorizations $Y^{(i)}=Q_YR_Y$ and $Z^{(i)}=Q_ZR_Z$ and an SVD $R_YR_Z^T=\hat Y\text{diag}(\sigma_1,\ldots,\sigma_\kappa)\hat Z^T$ where $\sigma_1,\ldots,\sigma_\kappa$ are the singular values in decreasing order. We can truncate to a rank-$\tilde\kappa$ matrix by dropping the terms corresponding to small singular values with a relative criterion $\sqrt{\sigma_{\tilde\kappa+1}^2+\cdots+\sigma_\kappa^2} \leq \epsilon_\text{rel}\sqrt{\sigma_1^2+\cdots+\sigma_\kappa^2}$ or an absolute one $\tilde\kappa=\max\{\tilde\kappa\mid\sigma_{\tilde\kappa}\geq\epsilon_\text{abs}\}$. In \textsc{Matlab} notation, the truncated matrix is $\tilde X^{(i)}=\tilde Y^{(i)}\tilde Z^{(i)T}$ with
\begin{equation}
    \tilde Y^{(i)} = Q_Y\hat Y(:,1:\tilde \kappa),\quad \tilde Z^{(i)} = Q_Z\hat Z(:,1:\tilde \kappa)\text{diag}(\sigma_1,\ldots,\sigma_{\tilde \kappa}).
\end{equation}

Low-rank approximation and truncation have been used for Krylov subspace methods \cite{BeOn15,KrTo11,LeEl17} and multigrid methods \cite{ElSu17}. More details can be found in these references. We will use examples of such solvers for linear systems arising in eigenvalue computations, as discussed in \cref{sec:diff,sec:stokes}.

\subsection{Normalization}
\label{sec:normalization}
In \cref{alg:sisi}, if $n_e=1$, the solution $v^{s,(i+1)}(\xi)$ is normalized so that $||u^{s,(i+1)}(\xi)||_2=1$ almost surely. With the superscripts omitted, assume $u(\xi)=\sum_{j=1}^{n_\xi}u_j\psi_j(\xi)$ is the normalized random vector constructed from $v(\xi)$. This expansion can be computed using sparse grid quadrature $\{\xi^{(q)},\eta^{(q)}\}_{q=1}^{n_q}$, where $\{\eta^{(q)}\}$ are the weights \cite{GeGr98}:
\begin{equation}
    \label{eq:normal_vec}
    u_j = \langle u(\xi)\psi_j(\xi)\rangle = \left\langle \frac{v(\xi)}{\| v(\xi)\|_2} \psi_j(\xi) \right\rangle \approx \sum_{q=1}^{n_q} \frac{v(\xi^{(q)})}{\| v(\xi^{(q)})\|_2} \psi_j(\xi^{(q)})\eta^{(q)}.
\end{equation}
Suppose the ``matricized'' version of the expansion coefficients of $v(\xi)$ is represented in low-rank form 
\begin{equation}
    \label{eq:lr_vu}
    V = [v_1,v_2,\ldots,v_{n_\xi}] = Y_vZ_v^T,
\end{equation}
where $Y_v\in\mathbb{R}^{n_x\times\kappa_v}$, $Z_v\in\mathbb{R}^{n_\xi\times\kappa_v}$. With $\Psi(\xi^{(q)})=[\psi_1(\xi^{(q)}),\psi_2(\xi^{(q)}),\ldots,\psi_{n_{\xi}}(\xi^{(q)})]^T$, we have
\begin{equation}
    \label{eq:yz_form}
    v(\xi^{(q)}) = \sum_{j=1}^{n_\xi} v_j\psi_j(\xi^{(q)}) = V\Psi(\xi^{(q)}) = Y_vZ_v^T\Psi(\xi^{(q)}).
\end{equation}
Let $U = [u_1,u_2,\ldots,u_{n_\xi}]$. Then \cref{eq:normal_vec} yields
\begin{equation}
    [U]_{:,j} = u_j =  \sum_{q=1}^{n_q} \frac{Y_vZ_v^T\Psi(\xi^{(q)})}{\| Y_vZ_v^T\Psi(\xi^{(q)})\|_2} \psi_j(\xi^{(q)})\eta^{(q)},
\end{equation}
and
\begin{equation}
    U  =  \sum_{q=1}^{n_q} \frac{Y_vZ_v^T\Psi(\xi^{(q)})}{\| Y_vZ_v^T\Psi(\xi^{(q)})\|_2} \Psi(\xi^{(q)})^T\eta^{(q)}.
\end{equation}
Thus, the matrix $U$ can be expressed as an outer product of two low-rank matrices $U=Y_uZ_u^T$ with
\begin{equation}
    Y_u = Y_v\in\mathbb{R}^{n_x\times\kappa_v}, \quad Z_u = \sum_{q=1}^{n_q} \frac{\Psi(\xi^{(q)})(\Psi(\xi^{(q)})^TZ_v)}{\| Y_v(Z_v^T\Psi(\xi^{(q)}))\|_2} \eta^{(q)}\in\mathbb{R}^{n_\xi\times\kappa_v}.
\end{equation}
This implies that the expansion coefficients of the normalized vector $u(\xi)$ can be written as a low-rank matrix with the same rank as the analogous matrix associated with $v(\xi)$. The cost of computing $Z_u$ is $O((n_x+n_\xi)n_q \kappa_v)$.

In the general case where more than one eigenvector is computed $(n_e>1)$, a stochastic version of the Gram-Schmidt process is applied to compute an orthonormal set $\{u^{s,(i+1)}(\xi)\}_{s=1}^{n_e}$ \cite{MeGh14,SoEl16}. With the superscript $(i+1)$ omitted, the process is based on the following calculation
\begin{equation}
    \label{eq:gram-schmidt}
    u^s(\xi) = v^s(\xi) - \sum_{t=1}^{s-1} \chi^{ts}(\xi) =  v^s(\xi) - \sum_{t=1}^{s-1} \frac{\langle v^s(\xi),u^t(\xi)\rangle_{\mathbb{R}^{n_x}}}{ \langle u^t(\xi),u^t(\xi)\rangle_{\mathbb{R}^{n_x}}} u^t(\xi)
\end{equation}
for $s=2,\ldots,n_e$. If the random vectors are represented as low-rank matrices as in \cref{eq:lr_vu}, the computational cost of the process can be reduced from $O(n_xn_\xi n_q)$ to $O((n_x+n_\xi)n_q\kappa)$ for some $\kappa\ll\min(n_x,n_\xi)$.

\subsection{Rayleigh quotient}
\label{sec:rayleigh}
The Rayleigh quotient in step 7 of \cref{alg:sisi} is computed (only once) after convergence of the inverse subspace iteration to find the eigenvalues. Given a normalized eigenvector $u(\xi)$ of problem \cref{eq:eigpb}, the computation of the stochastic Rayleigh quotient
\begin{equation}
    \label{eq:rayleigh}
    \lambda(\xi)=u(\xi)^TA(\xi)u(\xi)
\end{equation}
involves two steps:
\begin{enumerate}
\item[(1)] Compute matrix-vector product $w(\xi)=A(\xi)u(\xi)$ where $w(\xi)=\sum_{k=1}^{n_\xi}w_k\psi_k(\xi)$ and $w_k = \langle Au\psi_k\rangle$. In Kronecker product form,
\begin{equation}
    \label{eq:rayleigh_1}
    \mathbf{w} = \sum_{l=0}^{m} (G_l\otimes A_l)\mathbf{u}.
\end{equation}
If $\mathbf{u}$ has low-rank representation $U=Y_uZ_u^T$, then
\begin{equation}
    W = \sum_{l=0}^{m} (A_lY_u)(Z_uG_l)^T.
\end{equation}
\item[(2)] Compute eigenvalue $\lambda(\xi)=u(\xi)^Tw(\xi)$ where $\lambda(\xi)=\sum_{r=1}^{n_\xi}\lambda_r\psi_r(\xi)$ and $\lambda_r=\langle u^Tw\psi_r\rangle$. Equivalently,
\begin{equation}
    \lambda_r =  \langle \tilde G_r, H\rangle_{\mathbb{R}^{n_\xi\times n_\xi}} = \sum_{j,k=1}^{n_\xi} [\tilde G_r]_{jk} H_{jk}
\end{equation}
where $H_{jk}=u_j^Tw_k$ and thus $H=U^TW=Z_u(Y_u^TY_w)Z_w^T$. The matrices $\{\tilde G_r\}_{r=1}^{n_\xi}$ are sparse with $[\tilde G_r]_{jk}=\langle\psi_r\psi_j\psi_k\rangle$. In fact, if the basis functions are written as products of univariate polynomials, i.e.,
\begin{equation}
    \psi_r(\xi) = \psi_{r_1}(\xi_1)\psi_{r_2}(\xi_2)\cdots\psi_{r_m}(\xi_m),
\end{equation}
then $[\tilde G_r]_{jk}$ is nonzero only if $|j_l-k_l|\leq r_l \leq j_l+k_l$ and $r_l+j_l+k_l$ is even for all $1\leq l\leq m$ \cite{ErUl10}. This observation greatly reduces the cost of assembling the matrices $\{\tilde G_r\}$. For example, if $m=11$, the degree of the gPC basis functions is $p\leq 3$, and $n_\xi = (m+p)!/(m!p!) = 364$, then with the above rule, a total of 31098 nonzero entries must be computed, instead of the much larger number $n_\xi^3=48228544$ if the sparsity of $\{\tilde G_r\}$ is not used. 
\end{enumerate}

\subsection{Convergence criterion}
To check convergence, we can look at the magnitude of the residual
\begin{equation}
    r^s(\xi) = A(\xi)u^s(\xi)-\lambda^s(\xi)u^s(\xi),\quad s=1,2,\ldots,n_e.
\end{equation}
Alternatively, without computing the Rayleigh quotient at each iteration, error assessment can be done using the relative difference of the gPC coefficients of two successive iterates, i.e.,
\begin{equation}
    \epsilon_{\Delta u}^{s,(i)} = \frac{1}{n_\xi} \sum_{k=1}^{n_\xi} \frac{\|u_k^{s,(i)}-u_k^{s,(i-1)}\|_2}{\|u_k^{s,(i-1)}\|_2}.
\end{equation}
However, in the case of clustered eigenvalues (that is, if two or more eigenvalues are close to each other), the convergence of the inverse subspace iteration for single eigenvectors will be slow. Instead, we look at the angle between the eigenspaces \cite{BjGo73} in two consecutive iterations
\begin{equation}
\theta^{(i)}(\xi) = \angle(\text{span}(u^{1,(i)}(\xi),\ldots,u^{n_e,(i)}(\xi)),\,\text{span}(u^{1,(i-1)}(\xi),\ldots,u^{n_e,(i-1)}(\xi))).
\end{equation}
The expected value $\mathbb{E}[\theta^{(i)}]$ is taken as error indicator and is also calculated using sparse grid quadrature
\begin{equation}
    \label{eq:err_angle}
    \epsilon_\theta^{(i)} = \mathbb{E}[\theta^{(i)}] \approx \sum_{q=1}^{n_q} \theta^{(i)}(\xi^{(q)})\eta^{(q)}.
\end{equation}
At each quadrature point, $\theta^{(i)}(\xi^{(q)})$ is evaluated by \textsc{Matlab} function \verb+subspace+ for the largest principle angle.

\section{Stochastic diffusion equation}
\label{sec:diff}

In this section we consider the following elliptic equation with Dirichlet boundary conditions
\begin{equation}
\left\{ \begin{aligned}
    -\nabla\cdot(a(x,\omega)\nabla u(x,\omega)) & =\lambda(\omega)u(x,\omega) && \text{in } \mathcal{D}\times\Omega\\
    u(x,\omega) & =0 && \text{on }\partial \mathcal{D}\times\Omega
\end{aligned} \right.
\end{equation}
where $\mathcal{D}$ is a two-dimensional spatial domain and $\Omega$ is a sample space. The uncertainty in the problem is introduced by the stochastic diffusion coefficient $a(x,\omega)$. Assume that $a(x,\omega)$ is bounded and strictly positive and admits a truncated KL expansion
\begin{equation}
    \label{eq:kl}
    a(x,\omega) = a_0(x) + \sum_{l=1}^m \sqrt{\beta_l}a_l(x)\xi_l(\omega),
\end{equation}
where $a_0(x)$ is the mean function, $(\beta_l,a_l(x))$ is the $l$th eigenpair of the covariance function, and $\{\xi_l\}$ are a collection of uncorrelated random variables. The weak form is to find $(u(x,\xi),\lambda(\xi))$ such that for any $v(x)\in H_0^1(\mathcal{D})$,
\begin{equation}
    \int_\mathcal{D} a(x,\xi)\nabla u(x,\xi)\cdot\nabla v(x)\text{d}x = \lambda(\xi)\int_\mathcal{D} u(x,\xi)v(x)\text{d}x 
\end{equation}
almost surely. 

Finite element discretization in the physical domain $\mathcal{D}$ with basis functions $\{\phi_i(x)\}$ gives 
\begin{equation}  
    \label{eq:eig_diff_gnr}
    {K}(\xi){u}(\xi) = \lambda(\xi){M}{u}(\xi)
\end{equation}
where ${K}(\xi)=\sum_{l=0}^m K_l\xi_l$, and 
\begin{equation}
\begin{aligned}
    {[K_l]}_{ij} & = \int_\mathcal{D} \sqrt{\beta_l}a_l(x)\nabla\phi_i(x)\cdot\nabla\phi_j(x)\text{d}x,\\
    {[{M}]}_{ij} & = \int_\mathcal{D} \phi_i(x)\phi_j(x)\text{d}x, \quad i,j=1,2,\ldots,n_x,
\end{aligned}	
\end{equation}
with $\beta_0=1$ and $\xi_0\equiv1$. The result is a generalized eigenvalue problem where the matrix $M$ on the right-hand side is deterministic. With the Cholesky factorization $M=LL^T$, \cref{eq:eig_diff_gnr} can be converted to standard form 
\begin{equation} 
    \label{eq:eig_diff_sd}
    A(\xi){w}(\xi) = \lambda(\xi){w}(\xi),
\end{equation}
where $A(\xi)=L^{-1}K(\xi)L^{-T}$, ${w}(\xi)=L^T{u}(\xi)$.

We use stochastic inverse subspace iteration to find $n_e$ minimal eigenvalues of \cref{eq:eig_diff_sd}. As discussed in \cref{sec:stoch_isi}, the linear systems to be solved in each iteration are in the form
\begin{equation}
    \label{eq:eig_diff_solve0}
    \sum_{l=0}^m (G_l\otimes (L^{-1}K_lL^{-T}))\mathbf{v}^{s,(i+1)} = \mathbf{u}^{s,(i)},\quad s=1,2,\ldots,n_e.
\end{equation}
Let $\mathbf{v}^{s,(i)} = (I\otimes L^T)\hat{\mathbf{v}}^{s,(i)}$. Then \cref{eq:eig_diff_solve0} is equivalent to
\begin{equation}
    \label{eq:eig_diff_solve}
    \sum_{l=0}^m (G_l\otimes K_l) \hat{\mathbf{v}}^{s,(i+1)} =  (I\otimes L)\mathbf{u}^{s,(i)}.
\end{equation}

\subsection{Low-rank multigrid}
\label{sec:lrmg}
We developed a low-rank geometric multigrid me\-thod in \cite{ElSu17} for solving linear systems with the same structure as \cref{eq:eig_diff_solve}. The complete algorithm for solving
\begin{equation}
    \mathscr{A}(X)= \sum_{l=0}^m K_l XG_l^T =F
\end{equation}
is given in \cref{alg:lrmg}. All the iterates are expressed in low-rank form, and truncation operations are used to compress the ranks of the iterates. $\mathcal{T}_\text{rel}$ and $\mathcal{T}_\text{abs}$ are truncation operators with a relative tolerance $\epsilon_\text{rel}$ and an absolute tolerance $\epsilon_\text{abs}$, respectively. In each iteration, one V-cycle is applied to the residual equation. On the coarse grids, coarse versions of $\{K_l\}$ are assembled while the matrices $\{G_l\}$ stay the same. The prolongation operator is $\mathscr{P}=I\otimes P$, where $P$ is the same prolongation matrix as in a standard geometric multigrid solver, and the restriction operator is $\mathscr{R}=I\otimes P^T$. The smoothing operator $\mathscr{S}$ is based on a stationary iteration, and is also a Kronecker product of two matrices. The grid transfer and smoothing operations do not affect the rank. For instance, for any matrix iterate in low-rank form $X^{(i)}=Y^{(i)}Z^{(i)T}$,
\begin{equation}
   \mathscr{P}(X^{(i)}) = (PY^{(i)})(IZ^{(i)})^T.
\end{equation}
On the coarsest grid ($h=h_0$), the system is solved with direct methods.

\begin{algorithm2e}
\caption{Low-rank multigrid method}
\label{alg:lrmg}
\SetFuncSty{textsc}
\SetNlSty{}{}{:}
\DontPrintSemicolon
\SetKw{Init}{initialization}
\SetKwProg{Func}{function}{}{end}
\SetKwFunction{Vcycle}{Vcycle}
\SetKwFunction{Smooth}{Smooth}
\Init: $i=0$, $R^{(0)}=F$ in low-rank format, $r_0=\lVert F\rVert_F$ \;
\While{$r>tol*r_0$ $\&$ $i\leq maxit$}{
	$C^{(i)}=$ \Vcycle{$A,0,R^{(i)}$} \;
	$\tilde X^{(i+1)}=X^{(i)}+C^{(i)},$\tabto{6cm}$X^{(i+1)}=\mathcal{T}_{\text{abs}}(\tilde X^{(i+1)})$ \;
	$\tilde R^{(i+1)} = F-\mathscr{A}(X^{(i+1)}),$\tabto{6cm}$R^{(i+1)}=\mathcal{T}_{\text{abs}}(\tilde R^{(i+1)})$ \;
	$r=\lVert R^{(i+1)}\rVert_F$, $i=i+1$ \;
}
\BlankLine
\Func{$X^h=$ \Vcycle{$A^h,X_0^h,F^h$}}{
	\uIf{$h==h_0$}{
		solve $\mathscr{A}^h(X^h)=F^h$ directly \;
	}
	\Else{
		$X^h=$ \Smooth{$A^h,X_0^h,F^h$} \;
		$\tilde R^h=F^h-\mathscr{A}^h(X^h),$\tabto{5.5cm}$R^h=\mathcal{T}_\text{rel}(\tilde R^h)$ \;
		${R}^{2h}=\mathscr{R}(R^h)$ \;
		${C}^{2h}=$ \Vcycle{${A}^{2h},0,R^{2h}$} \;
		$X^{h}=X^h+\mathscr{P}(C^{2h})$ \;
		$X^h=$ \Smooth{$A^h,X^h,F^h$} \;
	}
}
\BlankLine
\Func{$X=$ \Smooth{$A,X,F$}}{
	\For{$\nu$ steps}{
		$\tilde X = X+\mathscr{S}(F-\mathscr{A}(X))$, \tabto{5.5cm}$X=\mathcal{T}_{\text{rel}}(\tilde X)$ \;
	}
}
\end{algorithm2e}

\subsection{Rayleigh-Ritz refinement}
\label{sec:rayleigh_ritz}
It is known that in the deterministic case with a constant diffusion coefficient, \cref{eq:eig_diff_gnr} typically has repeated eigenvalues \cite{ElSi14}, for example, $\lambda^2=\lambda^3$. The parametrized versions of these eigenvalues in the stochastic problem will be close to each other. In the deterministic setting, Rayleigh-Ritz refinement is used to accelerate the convergence of subspace iteration when some eigenvalues have nearly equal modulus and the convergence to individual eigenvectors is slow \cite{St69, St01}. Assume that a Hermitian matrix $S$ has eigendecomposition
\begin{equation}
    S = V\Lambda V^T = V_1\Lambda_1 V_1^T + V_2\Lambda_2 V_2^T
\end{equation}
where $\Lambda=\text{diag}(\lambda^1,\lambda^2,\ldots,\lambda^{n_x})$ with eigenvalues in increasing order and $V=[V_1,V_2]$ is orthogonal. Let the column space of $Q$ be a good approximation to that of $V_1$. Such an approximation is obtained from the inverse subspace iteration. The Rayleigh-Ritz procedure computes
\begin{itemize}
\item[(1)] Rayleigh quotient $T=Q^TSQ$, and
\item [(2)] eigendecomposition $T = W\Sigma W^T$.
\end{itemize}
Then $\Sigma$ and $QW$ represent good approximations to $\Lambda_1$ and $V_1$. 

The stochastic inverse subspace iteration algorithm produces solutions $\{u_\text{SG}^s(\xi)\}$ expressed as gPC expansions as in \cref{eq:gpc_expand} and sample eigenvectors are easily computed. The sample eigenvalues are generated from the stochastic Rayleigh quotient \cref{eq:rayleigh}. However, in the case of poorly separated eigenvalues, the sample solutions obtained this way are not accurate enough. Experimental results that demonstrate this are given in \cref{sec:num_diff}, see \cref{table:accuracy_norr}. Instead, we use a version of the Rayleigh-Ritz procedure to generate sample eigenvalues and eigenvectors with more accuracy. Specifically, a parametrized Rayleigh quotient $T(\xi)$ is computed using the approach of \cref{sec:rayleigh}, with
\begin{equation}
    [T]_{st}(\xi) = u_\text{SG}^{s}(\xi)^TA(\xi)u_\text{SG}^{t}(\xi),\quad s,t=1,2,\ldots,n_e.
\end{equation}
Then one can sample the matrix $T$, and for each realization $\xi^{(r)}$, solve a small ($n_e\times n_e$) deterministic eigenvalue problem $T(\xi^{(r)}) = W(\xi^{(r)})\Sigma(\xi^{(r)}) W(\xi^{(r)})^T$ to get better approximations for the minimal eigenvalues and corresponding eigenvectors:
\begin{equation}
\label{eq:rr_sample}
\begin{aligned}
\tilde\lambda^s_\text{SG}(\xi^{(r)}) & = [\Sigma(\xi^{(r)})]_{ss},\\
\tilde u^s_\text{SG}(\xi^{(r)}) & = [u^1_\text{SG}(\xi^{(r)}),u^2_\text{SG}(\xi^{(r)}),\ldots,u^{n_e}_\text{SG}(\xi^{(r)})] [W(\xi^{(r)})]_{:,s}.
\end{aligned}
\end{equation}
The effectiveness of this procedure will also be demonstrated in \cref{sec:num_diff}, see \cref{table:accuracy_rr}.

\subsection{Numerical experiments}
\label{sec:num_diff}
Consider a two-dimensional domain $\mathcal{D}=[-1,1]^2$. Let the spatial discretization consist of piecewise bilinear basis functions on a uniform square mesh.  The number of spatial degrees of freedom is $n_x = (2/h-1)^2$ where $h$ is the mesh size. Define the grid level $n_c$ such that $2/h=2^{n_c}$. In the KL expansion \cref{eq:kl}, we use an exponential covariance function
\begin{equation}
    r(x,y) = \sigma^2\text{exp}\left( -\frac{1}{b}\lVert x-y\rVert_1\right)
\end{equation}
and $(\beta_l,a_l(x))$ is the $l$th eigenpair of $r(x,y)$. The correlation length $b$ affects the decay of the eigenvalues $\{\beta_l\}$. The number of random variables $m$ is chosen so that $(\sum_{l=1}^m\beta_l)/(\sum_{l=1}^\infty\beta_l)\geq95\%$. Take the standard deviation $\sigma=0.01$, the mean function $a_0(x)\equiv1.0$, and $\{\xi_l\}$ to be independent and uniformly distributed on $[-\sqrt{3},\sqrt{3}]^m$. Legendre polynomials are used for gPC basis functions, whose total degree does not exceed $p=3$. The number of gPC basis functions is $n_\xi=(m+p)!/(m!p!)$. For the quadrature rule in \cref{sec:normalization}, we use a Smolyak sparse grid with Gauss Legendre quadrature points and grid level 4. For $m=11$, the number of sparse grid points is 2096.

\begin{figure}[hbtp]
    \label{fig:diff}
    \centering
    \subfloat[]{\label{fig:diff_ev}\includegraphics[width=0.49\textwidth]{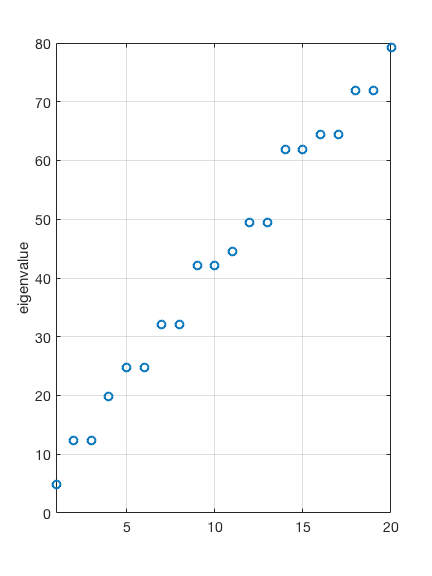}}
    \subfloat[]{\label{fig:diff_inexact}\includegraphics[width=0.49\textwidth]{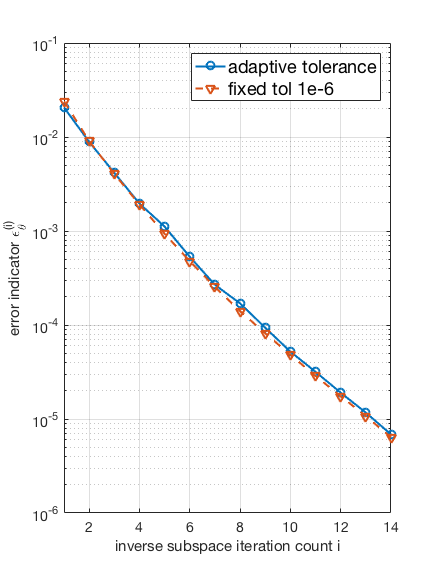}}
    \caption{(a): Smallest 20 eigenvalues of the mean problem. (b): Reduction of the error indicator $\epsilon_\theta^{(i)}$ for an adaptive multigrid tolerance \cref{eq:mg_tol_inexact} and a fixed tolerance $tol_\text{mg}=10^{-6}$. $n_c=6$, $m=11$.}
\end{figure}

We apply low-rank stochastic inverse subspace iteration to compute three minimal eigenvalues $(n_e=3)$ and corresponding eigenvectors for  \cref{eq:eig_diff_gnr}. The smallest 20 eigenvalues for the mean problem $K_0u=\lambda Mu$ are plotted in \cref{fig:diff_ev}. For the stochastic problem, the three smallest eigenvalues consist of one isolated smallest eigenvalue $\lambda^1(\xi)$ and (as mentioned in the previous subsection) two eigenvalues $\lambda^2(\xi)$ and $\lambda^3(\xi)$ that have nearly equal modulus. For the inverse subspace iteration, we take $\epsilon_\theta^{(i)}$ in \cref{eq:err_angle} as error indicator and use a stopping criterion $\epsilon_\theta^{(i)}\leq tol_\text{isi}=10^{-5}$. The low-rank multigrid method of \cref{sec:lrmg} is used to solve the system \cref{eq:eig_diff_solve}, where damped Jacobi iteration is employed for the smoothing opeator $\mathscr{S} = \omega_s \text{diag}(\mathscr{A})^{-1} = \omega_s(I\otimes K_0^{-1})$ with weight $\omega_s=2/3$. Two smoothing steps are applied ($\nu=2$). We also use the idea of inexact inverse iteration methods \cite{GoYe00,Lai97} so that in the first few steps of subspace iteration, the systems \cref{eq:ii_solve} are solved with milder error tolerances than in later steps. Specifically, we set the multigrid tolerance as
\begin{equation}
    \label{eq:mg_tol_inexact}
    tol_\text{mg}^{(i)} = \max\{ \min\{ 10^{-2}*\epsilon_\theta^{(i-1)},10^{-3}\}, 10^{-6}\},
\end{equation}
and truncation tolerances $\epsilon_\text{abs}^{(i)}=10^{-2}*tol_\text{mg}^{(i)}$, $\epsilon_\text{rel}=10^{-2}$ \cite{ElSu17}. This is shown to be useful in reducing the computational costs while not affecting the convergence of the subspace iteration algorithm (see \cref{fig:diff_inexact}). 

\Cref{table:rank} shows the ranks of the multigrid solutions in each iteration. It indicates that all the systems solved have low-rank approximate solutions ($n_x=3969$, $n_\xi=364$). With the inexact solve, the solutions have much smaller ranks in the first few iterations. In the last row of \cref{table:rank} are the numbers of multigrid steps $it_\text{mg}$ required to solve \cref{eq:eig_diff_solve} for $s=1$; similar numbers of multigrid steps are required for $s=2,3$. In addition, in \cref{alg:sisi} an absolute truncation operator with $\epsilon_\text{abs}=10^{-8}$ is applied after the computations in \cref{eq:gram-schmidt,eq:rayleigh_1} (both require addition of quantities represented as low-rank matrices in implementation) to compress the iterate ranks. Rayleigh-Ritz refinement discussed in \cref{sec:rayleigh_ritz} is used to obtain good approximations to individual sample eigenpairs.
\begin{table}[htbp!]
\caption{Iterate ranks after the multigrid solve and numbers of multigrid steps required in the inverse subspace iteration algorithm. $n_c=6$, $m=11$.}
\label{table:rank}
\centering
\begin{tabular}{c|c|*{14}{c}}
\hline
\multicolumn{2}{c|}{$(i)$} & 1 & 2 & 3 & 4 & 5 & 6 & 7 & 8 & 9 & 10 & 11 & 12 & 13 & 14\\ \hline
{\multirow{3}{*}{Rank}} & $u^1$ & 11 & 12 & 12 & 12 & 22 & 29 & 39 & 44 & 44 & 49 & 49 & 49 & 49 & 49\\ 
& $u^2$ & 11 & 12 & 13 & 14 & 23 & 28 & 43 & 49 & 52 & 54 & 54 & 54 & 54 & 54\\ 
& $u^3$ & 11 & 12 & 13 & 14 & 20 & 27 & 42 & 46 & 50 & 53 & 53 & 52 & 52 & 52\\ \hline
\multicolumn{2}{c|}{$it_\text{mg}$} & 3 & 4 & 4 & 5 & 5 & 5 & 6 & 6 & 6 & 7 & 7 & 7 & 7 & 7 \\ \hline
\end{tabular}
\end{table}

To show the accuracy of the low-rank stochastic Galerkin solutions, we compare them with reference solutions from Monte Carlo simulation. The stochastic Galerkin method produces a surrogate stochastic solution expressed with gPC basis functions that can be easily sampled. The Monte Carlo solutions are computed by the \verb+eigs+ function from \textsc{Matlab}, which uses the implicitly restarted Arnoldi method to compute several minimal eigenvalues \cite{So92}. For both methods, we use the same sample values $\{\xi^{(r)}\}$ of the random variables to generate sample eigenvalues and eigenvectors. Define the relative errors
\begin{equation}
\label{eq:relative_error}
\begin{aligned}
    \epsilon_{\lambda^s} &= \frac{1}{n_r} \sum_{r=1}^{n_r} \frac{| \lambda_\text{SG}^s(\xi^{(r)})-\lambda_\text{MC}^s(\xi^{(r)})|}{|\lambda_\text{MC}^s(\xi^{(r)})|}, \\
    \epsilon_{u^s} &= \frac{1}{n_r} \sum_{r=1}^{n_r} \frac{\| u_\text{SG}^s(\xi^{(r)})-u_\text{MC}^s(\xi^{(r)})\|_2}{\|u_\text{MC}^s(\xi^{(r)})\|_2},
\end{aligned}
\end{equation}
where $\lambda_\text{SG}^s$ and $u_\text{SG}^s$ denote the stochastic Galerkin sample solutions (they are replaced by $\tilde \lambda_\text{SG}^s$ and $\tilde u_\text{SG}^s$ in \cref{eq:rr_sample} if Rayleigh-Ritz refinement is used), $\lambda_\text{MC}^s$ and $u_\text{MC}^s$ are the Monte Carlo solutions, $n_r$ is the sample size, and $s=1,2,\ldots,n_e$. We use a sample size $n_r=10000$.  

We examine the accuracy for the three smallest eigenvalues obtained from inverse subspace iteration when they are computed both with and without Rayleigh-Ritz refinement. \cref{table:accuracy_norr} shows the results (for one spatial mesh size) when Rayleigh-Ritz refinement is not used. It can be seen that (the poorly separated) eigenvalues $\lambda^2$ and $\lambda^3$ are significantly less accurate than $\lambda^1$, and that the eigenvectors $u^2$ and $u^3$ are highly inaccurate. In contrast, \cref{table:accuracy_rr} (with results for three mesh sizes) demonstrates dramatically improved accuracy when refinement is done.  In all cases, convergence takes 14 iterations.

\begin{table}[thbp!]
\caption{Relative differences between low-rank stochastic Galerkin solutions (without Rayleigh-Ritz refinement) and Monte Carlo solutions. $nc=6$, $m=11$.}
\label{table:accuracy_norr}
\centering
\begin{tabular}{c|c||c|c}
\hline
$\epsilon_{\lambda^1}$ & $4.8393\times 10^{-10}$ & $\epsilon_{u^1}$ & $2.1494\times 10^{-7}$ \\ \hline
$\epsilon_{\lambda^2}$ & $1.4036\times 10^{-4}$ & $\epsilon_{u^2}$ & $3.2531\times 10^{-1}$ \\ \hline
$\epsilon_{\lambda^3}$ & $1.4021\times 10^{-4}$ & $\epsilon_{u^3}$ & $3.2530\times 10^{-1}$ \\
\hline
\end{tabular}
\end{table}

\begin{table}[htbp!]
\caption{Relative differences between low-rank stochastic Galerkin solutions (with Rayleigh-Ritz refinement) and Monte Carlo solutions. $m=11$.}
\label{table:accuracy_rr}
\centering
\begin{tabular}{c|c|c|c}
\hline
$n_c$ & 6 & 7 & 8 \\ \hline\hline
$\epsilon_{\lambda^1}$ & $4.8399\times 10^{-10}$ & $4.8443\times 10^{-10}$ & $4.8443\times 10^{-10}$ \\ \hline
$\epsilon_{\lambda^2}$ & $1.8641\times 10^{-9}$ & $1.8457\times 10^{-9}$ & $1.7780\times 10^{-9} $\\ \hline
$\epsilon_{\lambda^3}$ &$ 1.7027\times 10^{-9}$ &  $1.6849\times 10^{-9}$ & $1.6599\times 10^{-9}$ \\ \hline\hline
$\epsilon_{u^1}$ & $1.1390\times 10^{-7}$ & $1.8688\times 10^{-7}$ & $3.8872\times 10^{-7}$ \\ \hline
$\epsilon_{u^2}$ & $5.8266\times 10^{-6}$ & $6.0491\times 10^{-6}$ & $6.4745\times 10^{-6} $\\ \hline
$\epsilon_{u^3}$ & $5.8641\times 10^{-6}$ & $6.0874\times 10^{-6}$ & $6.5181\times 10^{-6}$ \\ \hline
\end{tabular}
\end{table}

The efficiency of the low-rank algorithm is demonstrated by comparison with (i) stochastic inverse subspace iteration without using low-rank approximation, which we call the full-rank stochastic Galerkin method, with the same tolerances $tol_\text{isi}$ and $tol_\text{mg}$, and (ii) the Monte Carlo method with a stopping tolerance $10^{-5}$ for the \verb+eigs+ function. The cost for computing stochastic Galerkin solutions consists of two parts, $t_\text{solve}$, the time required by inverse subspace iteration to compute the parametrized solution, and $t_\text{sample}$, the time to produce the sample solutions. (In the parlance of reduced basis methods \cite{VePa05}, the first part can be viewed as an offline computation.) As shown in \cref{table:time_nc}, once the parametrized solution is available, it is inexpensive to generate the sample solutions. The low-rank approximation greatly reduces both $t_\text{solve}$ and $t_\text{sample}$ of the stochastic Galerkin approach, especially as the mesh size gets refined. The total time required by low-rank stochastic Galerkin is much less than that for the Monte Carlo method, whereas the full-rank counterpart can be more expensive than Monte Carlo.

\begin{table}[htbp!]
\caption{Time comparison (in seconds) between stochastic Galerkin method and Monte Carlo simulation for various $n_c$. $b=4.0$, $m=11$, $n_\xi=364$.}
\label{table:time_nc}
\centering
\begin{tabular}{c|c|c|c|c}
\hline
\multicolumn{2}{c|}{$n_c$}  & 6 & 7 & 8 \\ \hline
\multicolumn{2}{c|}{$n_x$}  & 3969 & 16129 & 65025 \\ \hline \hline
\multirow{2}{*}{low-rank SG} & $t_\text{solve}$ & 262.50 & 852.86 & 3261.03 \\
& $t_\text{sample}$ &5.49 &16.07 & 78.94 \\ \hline
\multirow{2}{*}{full-rank SG} & $t_\text{solve}$ & 464.78 & 2230.01 & 20018.75 \\
& $t_\text{sample}$ &25.22 &104.46 & 422.36 \\ \hline
\multicolumn{2}{c|}{MC} & 500.73 & 2154.56 & 11563.40 \\\hline
\end{tabular}
\end{table}

\Cref{table:time_m} shows the performance of the stochastic Galerkin approach for various $n_\xi$, the number of degrees of freedom in the stochastic part. As expected, the Monte Carlo method is basically unaffected by the number of random variables in the KL expansion, whereas the cost of the stochastic Galerkin method increases as the number of parameters $m$ increases. The advantage of the stochastic Galerkin approach is clearer in the cases where $m$ is moderate. As the spatial mesh is refined (e.g., $n_c=8$), the low-rank stochastic Galerkin method becomes more efficient compared to Monte Carlo in all cases. The effectiveness of the low-rank algorithm is obvious for $m=16$ where the full-rank stochastic Galerkin method becomes too expensive or requires too much memory.

\begin{table}[htbp!]
\caption{Time comparison (in seconds) between stochastic Galerkin method and Monte Carlo simulation for various $m$.}
\label{table:time_m}
\centering
\subfloat[$n_c=7, n_x=16129$]{
\begin{tabular}{c|c|c|c|c}
\hline
\multicolumn{2}{c|}{$m$($b$)} & 8(5.0) & 11(4.0) & 16(3.0) \\ \hline
\multicolumn{2}{c|}{$n_\xi$} & 165 & 364 & 969 \\ \hline \hline
\multirow{2}{*}{low-rank SG} & $t_\text{solve}$ & 345.54 & 852.86 & 2627.54 \\ 
& $t_\text{sample}$ & 12.88 &16.07 & 22.80 \\ \hline
\multirow{2}{*}{full-rank SG} & $t_\text{solve}$ & 662.66 & 2230.01 & 13103.24 \\ 
& $t_\text{sample}$ & 46.02 &104.46 & 253.55 \\ \hline
\multicolumn{2}{c|}{MC} & 2069.57 & 2154.56 & 2394.18 \\\hline
\end{tabular}}

\subfloat[$n_c=8, n_x=65025$]{
\begin{tabular}{c|c|c|c|c}
\hline
\multicolumn{2}{c|}{$m$($b$)}  & 8(5.0) & 11(4.0) & 16(3.0) \\ \hline
\multicolumn{2}{c|}{$n_\xi$} & 165 & 364 & 969 \\ \hline \hline
\multirow{2}{*}{low-rank SG} & $t_\text{solve}$ & 1553.11 & 3261.03 & 8474.34 \\
& $t_\text{sample}$ &56.88 &78.94 & 100.32 \\ \hline
\multirow{2}{*}{full-rank SG} & $t_\text{solve}$ & 4725.65 &20018.75 &  out of \\
& $t_\text{sample}$ & 211.88 & 422.36 & memory \\ \hline
\multicolumn{2}{c|}{MC} & 11254.36 & 11563.40 & 12047.61 \\\hline
\end{tabular}}
\end{table}

\section{Stochastic Stokes equation}
\label{sec:stokes}

The second example of a stochastic eigenvalue problem that we consider is used to estimate the inf-sup stability constant associated with a discrete stochastic Stokes problem. Consider the following stochastic incompressible Stokes equation in a two-dimensional domain
\begin{equation}
\left\{\begin{aligned}
    -\nabla\cdot(a(x,\omega)\nabla\vec{u}(x,\omega))+\nabla p(x,\omega) &=\vec{0} && \text{ in } \mathcal{D}\times\Omega\\
    \nabla\cdot\vec{u}(x,\omega) &= 0 && \text{ in } \mathcal{D}\times\Omega
\end{aligned}\right.
\end{equation}
with a Dirichlet inflow boundary condition $\vec{u}(x,\omega) = \vec{u}_D(x)$ on $\partial\mathcal{D}_D\times\Omega$ and a Neumann outflow boundary condition $a(x,\omega)\nabla\vec{u}(x,\omega)\cdot\vec{n} - p(x,\omega)\vec{n} = \vec{0}$ on $\partial\mathcal{D}_N\times\Omega$. Such problems and more general stochastic Navier--Stokes equations have been studied in  \cite{PoSi12,SoEl16b}. As in the diffusion problem, we assume that the stochastic viscosity $a(x,\omega)$ is represented by a truncated KL expansion \cref{eq:kl} with random variables $\{\xi_l\}_{l=1}^m$. The weak formulation of the problem is: find $\vec{u}(x,\xi)$ and $p(x,\xi)$ satisfying
\begin{equation}
\left\{\begin{aligned}
    \int_\mathcal{D} a(x,\xi) \nabla\vec{u}(x,\xi):\nabla\vec{v}(x) - p(x,\xi)\nabla\cdot\vec{v}(x)\,\text{d}x &= 0\\
    \int_\mathcal{D} q(x)\nabla\cdot\vec{u}(x,\xi)\,\text{d}x & = 0
\end{aligned} \right.
\end{equation}
almost surely for any $\vec{v}(x)\in H^1_0(\mathcal{D})^2$ (zero boundary conditions on $\partial\mathcal{D}_D$) and $q(x)\in L^2(\mathcal{D})$. Here $\nabla\vec{u}:\nabla\vec{v}$ is a componentwise scalar product ($\nabla u_{x_1}\cdot\nabla v_{x_1}+\nabla u_{x_2}\cdot\nabla v_{x_2}$ for two-dimensional $(u_{x_1},u_{x_2})$). Finite element discretization with basis functions $\{\vec{\phi}_i(x)\}$ for the velocity field and $\{\varphi_k(x)\}$ for the pressure field results in a linear system in the form 
\begin{equation}
    \begin{pmatrix}
    K(\xi) & B^T\\
    B & 0\\
    \end{pmatrix}\begin{pmatrix}
    \vec{u}(\xi)\\
    p(\xi)
    \end{pmatrix} = \begin{pmatrix}
    f\\g
    \end{pmatrix},
\end{equation}
where $K(\xi)=\sum_{l=0}^m K_l\xi_l$, and
\begin{equation}
\begin{aligned}
    {[K_l]}_{ij} &= \int_\mathcal{D} \sqrt{\beta_l}a_l(x) \nabla\vec{\phi}_i(x):\nabla\vec{\phi}_j(x)\text{d}x,\\
    {[B]}_{kj} &= -\int_\mathcal{D} \varphi_k(x)\nabla\cdot\vec{\phi}_j(x)\text{d}x,
\end{aligned}
\end{equation}
for $i,j=1,2,\ldots,n_u$ and $k=1,2,\ldots,n_p$. The Dirichlet boundary condition is incorporated in the right-hand side. 

We are interested in the parametrized inf-sup stability constant $\gamma(\xi)$ for the discrete problem. Evaluation of the inf-sup constant for various parameter values plays an important role for \textit{a posteriori} error estimation for reduced basis methods \cite{Cu05,VePa05}. For this, we exploit the fact that  $\gamma(\xi)$ has an algebraic interpretation \cite{ElSi14}
\begin{equation}
    \label{eq:infsup}
    \gamma^2(\xi) = \min_{q(\xi)\neq 0}\frac{\langle BK(\xi)^{-1}B^Tq(\xi),q(\xi)\rangle_{\mathbb{R}^{n_p}}}{\langle Mq(\xi),q(\xi)\rangle_{\mathbb{R}^{n_p}}}
\end{equation}
where $M$ is the mass matrix with ${[M]}_{ij} = \int_\mathcal{D} \varphi_i(x)\varphi_j(x)$, $i,j=1,2,\ldots,n_p$. Thus, finding $\gamma(\xi)$ is equivalent to finding the smallest eigenvalue of the generalized eigenvalue problem
\begin{equation}
    \label{eq:stokes_eig}
    BK(\xi)^{-1}B^Tq(\xi) = \lambda(\xi)Mq(\xi)
\end{equation}
associated with the stochastic pressure Schur complement $BK(\xi)^{-1}B^T$. This can be written in standard form as
\begin{equation}
   \label{eq:stokes_eig_sd}
   L^{-1}BK(\xi)^{-1}B^TL^{-T}w(\xi) = \lambda(\xi)w(\xi)
\end{equation}
where $M=LL^T$ is a Cholesky factorization, and $w(\xi)=L^Tq(\xi)$. 

The eigenvalue problem \cref{eq:stokes_eig_sd} does not have exactly the same form as \cref{eq:eigpb}, since it involves the inverse of $K(\xi)$. If we use the stochastic inverse iteration algorithm to compute the minimal eigenvalue of \cref{eq:stokes_eig_sd}, then each iteration requires solving
\begin{equation}
    \label{eq:stokes_ii_solve}
    \langle L^{-1}BK^{-1}B^TL^{-T} v^{(i+1)}\psi_k\rangle = \langle u^{(i)}\psi_k\rangle,\quad k=1,2,\ldots,n_\xi,
\end{equation}
for $v^{(i+1)}(\xi)$. We can reformulate \cref{eq:stokes_ii_solve} to take advantage of the Kronecker product structure and low-rank solvers. Let $z(\xi)=-K(\xi)^{-1}B^TL^{-T}{v}^{(i+1)}(\xi)$ and let $\hat{v}^{(i+1)}(\xi)=L^{-T}{v}^{(i+1)}(\xi)$. Then \cref{eq:stokes_ii_solve} is equivalent to the coupled system
\begin{equation}
    \label{eq:stokes_aux1}
    \langle (Kz+B^T\hat{v}^{(i+1)})\psi_k\rangle=0,\quad \langle Bz\psi_k\rangle = \langle -Lu^{(i)}\psi_k\rangle,\quad k=1,2,\ldots,n_\xi.
\end{equation}
As discussed in \cref{sec:stoch_isi}, the random vectors are expressed as gPC expansions. Thus, \cref{eq:stokes_aux1} can be written in Kronecker product form as a discrete Stokes system for coefficient vectors $\mathbf{z}$, $\hat{\mathbf{v}}^{(i+1)}$,
\begin{equation}
    \label{eq:stokes_galerkin}
    \begin{pmatrix}
    \sum_{l=0}^m ( G_l\otimes K_l )& I\otimes B^T\\
    I\otimes B & 0 \end{pmatrix}
    \begin{pmatrix}
    \mathbf{z} \\ 
    \hat{\mathbf{v}}^{(i+1)} \end{pmatrix}
    = \begin{pmatrix}
    \mathbf{0}\\
    -(I\otimes L)\mathbf{u}^{(i)} \end{pmatrix},
\end{equation}
and ${\mathbf{v}}^{(i+1)}=(I\otimes L^T)\hat{\mathbf{v}}^{(i+1)}$.

In addition, for the eigenvalue problem \cref{eq:stokes_eig_sd}, computing the Rayleigh quotient \cref{eq:rayleigh} requires solving a linear system. In the first step of \cref{eq:rayleigh}, for the matrix-vector product, one needs to compute $w(\xi)=K(\xi)^{-1}\hat{u}(\xi)$, where $\hat{u}(\xi)=B^TL^{-T}u(\xi)$. For the weak formulation, this corresponds to solving a linear system
\begin{equation}
    \label{eq:stokes_rayleigh}
     \left(\sum_{l=0}^m G_l\otimes K_l\right) \mathbf{w} = \hat{\mathbf{u}}.
\end{equation}

\subsection{Low-rank MINRES}
We discuss a low-rank iterative solver for \cref{eq:stokes_galerkin}. The system is symmetric but indefinite, with a positive-definite $(1,1)$ block. A low-rank preconditioned MINRES method for solving $\mathscr{A}(X)=F$ is used and described in \cref{alg:minres}. The precondtioner is block-diagonal
\begin{equation}
    \mathscr{M} = \begin{pmatrix}
    \mathscr{M}_{11} & 0\\
    0  &   \mathscr{M}_{22} \end{pmatrix}.
\end{equation}
We use an approximate mean-based preconditioner \cite{PoEl09} for the $(1,1)$ block: $\mathscr{M}_{11} = G_0\otimes \hat K_0$ = $I\otimes \hat K_0$. Here, $\hat K_0^{-1}$ is defined by approximation of the action of $K_0^{-1}$, using one V-cycle of an algebraic multigrid method (AMG) \cite{SiSi11}. For the $(2,2)$ block, we take $\mathscr{M}_{22}=I\otimes \text{diag}(M)$. As in the multigrid method, all the quantities are in low-rank format, and truncation operations are applied to compress matrix ranks. \cref{alg:minres} requires the computation of inner products of two low-rank matrices $\langle X_1,X_2\rangle_{\mathbb{R}^{n_x\times n_\xi}}$. Let $X_1=Y_1Z_1^T$, $X_2=Y_2Z_2^T$ with $Y_1\in\mathbb{R}^{n_x\times\kappa_1}$, $Z_1\in\mathbb{R}^{n_\xi\times\kappa_1}$, $Y_2\in\mathbb{R}^{n_x\times\kappa_2}$, $Z_2\in\mathbb{R}^{n_\xi\times\kappa_2}$. Then the inner product can be computed with a cost of $O((n_x+n_\xi+1)\kappa_1\kappa_2)$ \cite{KrTo11}:
\begin{equation}
    \langle X_1,X_2\rangle = \text{trace}(X_1^TX_2)
    = \text{trace}(Z_1Y_1^TY_2Z_2^T)
    = \text{trace}((Z_2^TZ_1)(Y_1^TY_2)).
\end{equation}

\begin{algorithm2e}
\caption{Low-rank preconditioned MINRES method}
\label{alg:minres}
\SetNlSty{}{}{:}
\DontPrintSemicolon
\SetKw{Init}{initialization}
\Init: $V^{(0)}=0$, $W^{(0)}=0$, $W^{(1)}=0$, $\gamma_0=0$. Choose $X^{(0)}$, compute $V^{(1)}=F-\mathscr{A}(X^{(0)})$. $P^{(1)}=\mathscr{M}^{-1}(V^{(1)})$, $\gamma_1=\sqrt{\langle P^{(1)},V^{(1)}\rangle}$. Set $\eta=\gamma_1$, $s_0=s_1=0$, and $c_0=c_1=1$. \;
\For{$j=1,2,\ldots$}{
    $P^{(j)} = P^{(j)}/\gamma_j$ \;
    $\tilde R^{(j)} =\mathscr{A}(P^{(j)})$, \tabto{7.6cm} $R^{(j)}=\mathcal{T}_\text{rel}(\tilde R^{(j)})$ \label{ln:4}\;
    $\delta_j = \langle R^{(j)},P^{(j)}\rangle$ \;
    $\tilde V^{(j+1)} = R^{(j)} - (\delta_j/\gamma_j)V^{(j)}- (\gamma_j/\gamma_{j-1})V^{(j-1)}$,\tabto{7.6cm} $V^{(j+1)}=\mathcal{T}_\text{rel}(\tilde V^{(j+1)})$ \label{ln:6}\;
    $P^{(j+1)} = \mathscr{M}^{-1}(V^{(j+1)})$ \;
    $\gamma_{j+1}=\sqrt{\langle P^{(j+1)},V^{(j+1)}\rangle}$ \;
    $\alpha_0 = c_j\delta_j - c_{j-1}s_j\gamma_j$ \;
    $\alpha_1 = \sqrt{\alpha_0^2 + \gamma_{j+1}^2}$ \;
    $\alpha_2 = s_j\delta_j + c_{j-1}c_j\gamma_j$ \;
    $\alpha_3 = s_{j-1}\gamma_j$ \;
    $c_{j+1} = \alpha_0/\alpha_1$, $s_{j+1}=\gamma_{j+1}/\alpha_1$ \;
    $\tilde W^{(j+1)} = (P^{(j)} - \alpha_3 W^{(j-1)} - \alpha_2 W^{(j)})/\alpha_1 $,\tabto{7.6cm} $W^{(j+1)}=\mathcal{T}_\text{rel}(\tilde W^{(j+1)})$ \label{ln:14}\;
    $\tilde X^{(j)}= X^{(j-1)} + c_{j+1}\eta W^{(j+1)}$,\tabto{7.6cm} $X^{(j)}=\mathcal{T}_\text{rel}(\tilde X^{(j)})$ \;
    $\eta = -s_{j+1}\eta$ \;
    Check convergence
}
\end{algorithm2e}

\subsection{Numerical experiments}
Consider a two-dimension channel flow on domain $\mathcal{D}=[-1,1]^2$ with uniform square meshes. Let $\partial\mathcal{D}_D=\{(x_1,x_2)\mid x_1=-1,\text{ or } x_2=1,\text{ or } x_2=-1\}$ and $\partial\mathcal{D}_N=\{(x_1,x_2)\mid x_1=1\}$. Define grid level $n_c$ so that $2/h=2^{n_c}$, where $h$ is the mesh size. We use the Taylor-Hood method for finite element discretization with biquadratic basis functions $\{\vec{\phi}_i(x)\}$ for the velocity field and bilinear basis functions $\{\varphi_k(x)\}$ for the pressure field. For the velocity field the basis functions are in the form $\left\{\left(\begin{smallmatrix}\phi_i(x)\\0\end{smallmatrix}\right),\left(\begin{smallmatrix}0\\\phi_i(x)\end{smallmatrix}\right)\right\}$, where $\{\phi_i(x)\}$ are scalar-value biquadratic basis functions. The number of degrees of freedom in the spatial discretization is $n_x=n_u+n_p$ where $n_p=(2^{n_c}+1)^2$. Assume the viscosity $a(x,\xi)$ has a KL expansion with the same specifications as in the diffusion problem. For the quadrature rule in \cref{sec:normalization}, we use a Smolyak sparse grid with Gauss Legendre quadrature points and grid level 4.

\begin{figure}[hbtp]
    \label{fig:stokes}
    \centering
    \subfloat[]{\label{fig:stokes_ev}\includegraphics[width=0.5\textwidth]{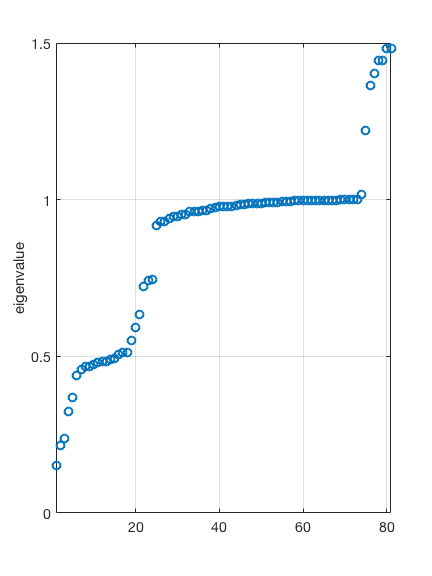}}
    \subfloat[]{\label{fig:stokes_minres}\includegraphics[width=0.5\textwidth]{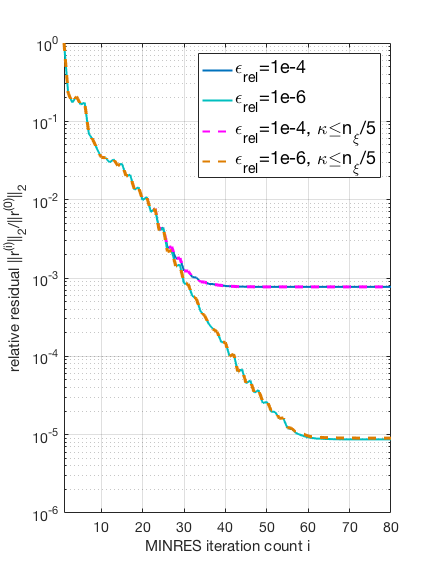}}
    \caption{(a): Eigenvalues of $BK_0^{-1}B^Tq=\lambda Mq$. $n_c=3$. (b): Reduction of the relative residual for the low-rank MINRES method with various truncation criteria. Solid lines: relative tolerance $\epsilon_\text{rel}$; dashed lines: relative tolerance $\epsilon_\text{rel}$ with rank $\kappa\leq n_\xi/5$. $n_c=4$, $m=11$.}
\end{figure}

We use stochastic inverse iteration algorithm to find the minimal eigenvalue of \cref{eq:stokes_eig}. The eigenvalues of $BK_0^{-1}B^Tq=\lambda Mq$ are plotted in \cref{fig:stokes_ev} with $n_c=3$. It shows that the minimal eigenvalue is isolated from the larger ones. For the inverse iteration, we take $\epsilon_\theta^{(i)}$ in \cref{eq:err_angle} as error indicator and use a stopping criterion $\epsilon_\theta^{(i)}\leq tol_\text{isi}=10^{-5}$. The error tolerance for the MINRES solver $tol_\text{minres}^{(i)}$ is set as in \cref{eq:mg_tol_inexact}. \Cref{fig:stokes_minres} shows the convergence of the low-rank MINRES method for different relative truncation tolerances $\epsilon_\text{rel}$. It indicates the accuracy that MINRES can achieve is related to $\epsilon_\text{rel}$. In the numerical experiments we use $\epsilon_\text{rel}^{(i)}=10^{-1}* tol_\text{minres}^{(i)}$. In addition, we have observed that in many cases the truncations in Lines \ref{ln:4}, \ref{ln:6}, \ref{ln:14} of \cref{alg:minres} produce relatively high ranks, which increases the computational cost. To handle this, we impose a bound on the ranks $\kappa$ of the outputs of these truncation operators such that $\kappa\leq n_\xi/5$ (in general $n_x\geq n_\xi$). It is shown in \Cref{fig:stokes_minres} that the convergence of low-rank MINRES is unaffected by this strategy. \Cref{table:stokes_rank} shows the ranks of the MINRES solution $\hat{\mathbf{v}}^{(i)} $ in \cref{eq:stokes_galerkin} and numbers of MINRES steps $it_\text{minres}$ required in each iteration. Note that the system being solved has size $n_x n_\xi\times n_x n_\xi$ whereas the matricized $\hat{V}^{(i)}\in\mathbb{R}^{n_p\times n_\xi}$. For $n_c=4$ and $m=11$, $n_p=289$, $n_x=2273$, $n_\xi=364$. The ranks of the solutions are no larger than 51. For the Rayleigh quotient, the system \cref{eq:stokes_rayleigh} is solved by a low-rank conjugate gradient method \cite{KrTo11} with a relative residual smaller than $10^{-8}$.

\begin{table}[htbp!]
\caption{Iterate ranks after the MINRES solve and numbers of MINRES steps required in the inverse iteration algorithm. $n_c=4$, $m=11$.}
\label{table:stokes_rank}
\centering
\setlength{\tabcolsep}{5pt}
\begin{tabular}{c|*{16}{c}}
\hline
$(i)$ & 1 & 2 & 3 & 4 & 5 & 6 & 7 & 8 & 9 & 10 & 11 & 12 & 13 & 14 & 15 & 16\\ \hline
Rank & 7 & 14 & 14 & 15 & 17 & 22 & 28 & 34 & 43 & 51 & 51 & 50 & 50 & 51 & 50 & 51\\ \hline 
$it_\text{minres}$ &  31 & 50 & 52 & 54 & 57 & 59 & 62 & 64 & 66 & 67 & 67 & 67 & 67 & 67 & 67 & 67 \\ \hline
\end{tabular}
\end{table}

As in the diffusion problem, we compare the results from the stochastic Galerkin approach with those from Monte Carlo simulations. Let $m=11$, $p=3$, $n_\xi=364$. We use a sample size $n_r=1000$. \Cref{table:stokes_eps} shows the accuracy of the stochastic Galerkin solutions where $\epsilon_{\lambda^1}$ and $\epsilon_{u^1}$ are defined in \cref{eq:relative_error} (no Rayleigh-Ritz procedure is used here). In all cases, convergence of the inverse iteration takes 16--18 steps. 

The efficiency of the low-rank algorithm is shown in \cref{table:stokes_time} by comparison with the full-rank stochastic Galerkin method and the Monte Carlo method (with a stopping tolerance $10^{-5}$). For the latter, the \verb+eigs+ function requires solving linear systems associated with $BK(\xi^{(r)})^{-1}B^T$ for each sample $\xi^{(r)}$. This can be achieved by computing the coefficient matrix and using a direct solver, or reformulating as a discrete Stokes problem which is solved by a preconditioned MINRES method with a relative residual $\leq 10^{-6}$. We use whichever is more efficient to make a fair comparison. From \cref{table:stokes_time} it is clear that the low-rank approximation enhances the efficiency of the stochastic Galerkin approach. The cost of computing stochastic Galerkin solutions increases more slowly than that of the Monte Carlo method as $n_c$ becomes larger. As the mesh gets refined, the stochastic Galerkin approach becomes more efficient compared to Monte Carlo. It should also be noted that we are using a relatively small sample size. For $n_r=10000$, the computing time of Monte Carlo simulations will be about 10 times larger; for the stochastic Galerkin approach only $t_\text{sample}$ is affected, and that extra cost will be negligible since $t_\text{sample}$ is so small. With a smaller $m=8$ in \cref{table:stokes_time2}, the stochastic Galerkin approach uses less time whereas Monte Carlo is basically unaffected.

\begin{table}[htbp!]
\caption{Relative difference between stochastic Galerkin solutions and Monte Carlo solutions. $b=4.0$, $m=11$, $n_\xi=364$.}
\label{table:stokes_eps}
\centering
\begin{tabular}{c|c|c|c}
\hline
$n_c$ & 4 & 5 & 6 \\ \hline
$\epsilon_{\lambda^1}$ & $5.8717\times 10^{-9}$ & $7.5615\times 10^{-9}$ & $9.2116\times 10^{-9}$ \\ \hline
$\epsilon_{u^1}$ & $4.9543\times 10^{-5}$ & $5.5996\times 10^{-5}$ & $5.7339\times 10^{-5}$ \\ \hline
\end{tabular}
\end{table}

\begin{table}[htbp!]
\caption{Time comparison (in seconds) between stochastic Galerkin method and Monte Carlo simulation for various $n_c$.}
\label{table:stokes_time}
\centering
\subfloat[$b=4.0$, $m=11$, $n_\xi=364$]{\label{table:stokes_time1}
\begin{tabular}{c|c|c|c|c}
\hline
\multicolumn{2}{c|}{$n_c$} & 4 & 5 & 6 \\ \hline
\multicolumn{2}{c|}{$n_p$} & 289 & 1089 & 4225 \\  \hline
\multicolumn{2}{c|}{$n_x$} & 2273 & 9153 & 36737 \\ \hline \hline
\multirow{2}{*}{low-rank SG} & $t_\text{solve}$ & 319.59 & 1245.72 & 5896.94 \\ 
& $t_\text{sample}$ &0.08 & 0.13 & 0.39 \\ \hline
\multirow{2}{*}{full-rank SG} & $t_\text{solve}$ & 417.87 & 1755.28 & 8179.09 \\ 
& $t_\text{sample}$ &0.09 &0.16 & 0.92 \\ \hline
\multicolumn{2}{c|}{MC} & 201.91 & 3847.30 & 17767.97 \\\hline
\end{tabular}}

\subfloat[$b=5.0$, $m=8$, $n_\xi=165$]{\label{table:stokes_time2}
\begin{tabular}{c|c|c|c|c}
\hline
\multicolumn{2}{c|}{$n_c$} & 4 & 5 & 6 \\ \hline
\multicolumn{2}{c|}{$n_p$} & 289 & 1089 & 4225 \\  \hline
\multicolumn{2}{c|}{$n_x$} & 2273 & 9153 & 36737 \\ \hline \hline
\multirow{2}{*}{low-rank SG} & $t_\text{solve}$ & 143.32 & 447.98 &  2509.00\\ 
& $t_\text{sample}$ & 0.06 & 0.06 & 0.10 \\ \hline
\multirow{2}{*}{full-rank SG} & $t_\text{solve}$ & 189.96 & 830.36 & 3810.33 \\ 
& $t_\text{sample}$ & 0.06 &0.07 &0.36 \\ \hline
\multicolumn{2}{c|}{MC} & 201.49 & 3847.90 & 17782.69 \\\hline
\end{tabular}}

\end{table}

\section{Summary} We studied low-rank solution methods for the stochastic eigenvalue problems. The stochastic Galerkin approach was used to compute surrogate approximations to the minimal eigenvalues and corresponding eigenvectors, which are stochastic functions with gPC expansions. We introduced low-rank approximations to enhance efficiency of the stochastic inverse subspace iteration algorithm. Two detailed benchmark problems, the stochastic diffusion problem, and an operator associated with a discrete stochastic Stokes equation, were considered for illustrating the effectiveness of the proposed low-rank algorithm. It was confirmed in the numerical experiments that the low-rank solution method produces accurate results with much less computing time, making the stochastic Galerkin method more competitive compared with the sample-based Monte Carlo approach.

\bibliographystyle{siamplain}
\bibliography{stoch_eig}

\end{document}